\documentclass{amsart}
\usepackage{amssymb}
\usepackage{graphicx}
\newtheorem{theorem}{Theorem}[section]
\newtheorem{lemma}[theorem]{Lemma}
\newtheorem{proposition}[theorem]{Proposition}
\newtheorem{corollary}[theorem]{Corollary}
\newtheorem*{acknowledgement}{Acknowledgement}
\theoremstyle{definition}
\newtheorem{definition}[theorem]{Definition}

\theoremstyle{remark}
\newtheorem{remark}[theorem]{Remark}

\numberwithin{equation}{section}

\begin{document}
\title[Invertible harmonic mappings]{Invertible harmonic mappings, beyond Kneser}
\author{Giovanni Alessandrini}
\address{Dipartimento di Matematica e Informatica, Universit\`a degli Studi di Trieste, Italy}

\email{alessang@units.it}
\thanks{The first author was supported  in part by MiUR, PRIN no. 2006014115}

\author{Vincenzo Nesi}
\address{Dipartimento di Matematica,
La Sapienza, Universit\`a di Roma, Italy}
\email{nesi@mat.uniroma1.it}
\thanks{The second author was supported in part by MiUR, PRIN
no. 2006017833.}


\begin{abstract}
We prove necessary and sufficient criteria of invertibility for
planar harmonic mappings which generalize a classical result of H.
Kneser, also known as the Rad\'{o}--Kneser--Choquet theorem.
\end{abstract}

\maketitle
\section{Introduction}
\label{intro} Let $B := \{(x,y)\in \mathbb R^2: x^2+y^2<1\}$ denote
the unit disk. Given a homeomorphism $\Phi$ from the unit circle
$\partial B$ onto a simple closed curve $\gamma\subseteq \mathbb
R^2$, let us consider the   solution $U\in C^2 (B;\mathbb R^2)\cap
C(\overline{B};\mathbb R^2)$ to the following Dirichlet problem
\begin{equation}\label{main}
\left\{
\begin{array}{lll}
\Delta U=0,&\hbox{in}&B,\\
U=\Phi,&\hbox{on}&\partial B.
\end{array}
\right.
\end{equation}
The basic question that we address in this paper is under which conditions on $\Phi$ we
have that $U$ is a homeomorphism of $\overline{B}$ onto $\overline{D}$, where $D$ denotes the bounded open,
simply connected  set for which $\partial D=\gamma$.

The fundamental benchmark for this issue is a classical theorem,
first conjectured by T. Rad\'o in 1926 \cite{rado1}, which was
proved immediately after by H. Kneser \cite{kneser}, and
subsequently rediscovered, with a different proof, by G. Choquet
\cite{Choquet}. Let us recall the result.
\begin{theorem}[H. Kneser]\label{Kneser.th}
If $D$ is convex, then $U$ is a homeomorphism of $\overline{B}$ onto $\overline{D}$.
\end{theorem}
We recall that this Theorem had a remarkable impact in the
development of the theory of minimal surfaces, see for instance
\cite{rado2}. Its   influence appears also in other areas of
mathematics, let us mention here homogenization and effective
properties of materials \cite{bmn, an:arch, an:cocv}, inverse
boundary value problems \cite{fv, adiaz, seo} and, quite recently,
variational problems for maps of finite distortion \cite{aio}. See
also, as general references, and for many interesting related
results, the book by Duren \cite{duren}  and the review article by
Bshouty and Hengartner \cite{bh}.

The amazing character of Kneser's Theorem stands in the simplicity
and elegance of the geometric condition on the target curve
$\gamma$. Let us emphasize here that this condition does {\em not}
involve the choice of the parametrization $\Phi$ of the curve
$\gamma$.

 In order to motivate the main result of this paper, Theorem \ref{main.th} below, we wish to stress that \emph{no} weaker condition on the shape of $D$ can
replace the assumption in Theorem \ref{Kneser.th}. In fact, the
following Theorem holds.
\begin{theorem}[G. Choquet]\label{parabola.th}
For every Jordan domain $D$ which is \emph{not} convex, there exists  a homeomorphism $\Phi:\partial B \to \partial D$
such that the solution $U$ to (\ref{main}) is \emph{not} a homeomorphism.
\end{theorem}

A proof for this Theorem is due to Choquet \cite[\S 3]{Choquet}. In
Section 6 we present a new proof aimed at having a more explicit
description of the homeomorphism $\Phi$.
 In the final part of this Introduction, when presenting the content of Section 6,
 we shall
illustrate the advantages of this new proof  with more details.

Theorem  \ref{parabola.th} shows that, given a non--convex domain
$D$ and its boundary $\gamma$, one can find some parameterization
 of the latter which give rise to a non--invertible solution to (\ref{main}). On the other hand, by the Riemann Mapping
 Theorem, see for
instance \cite[Theorem 3.4]{palka},
 for any such $\gamma$ one can also find other parameterizations  for which the corresponding solution to (\ref{main}) is a homeomorphism and,
 in fact, a conformal mapping. Thus the question arises, for a given simply connected target
domain $D$, possibly non--convex, of how to characterize {\em all}
the parameterizations which give rise to an invertible solution to
(\ref{main}).

Our main result is a complete answer to this question for those
parameterizations $\Phi$ which are smooth enough so that the
corresponding solution  to    (\ref{main}) belongs to
$C^1(\overline{B};\mathbb R^2)$.

\begin{theorem}\label{main.th}
Let $\Phi:\partial B\to \gamma \subset \mathbb R^2$ be an orientation preserving diffeomorphism of class $C^1$ onto a simple closed curve $\gamma$.
Let $D$ be the bounded domain such that $\partial D = \gamma$.
 Let $U\in C^2 (B;\mathbb R^2)\cap C(\overline{B};\mathbb R^2)$ be  the   solution to   (\ref{main}) and assume, in addition, that
 $U\in C^1(\overline{B};\mathbb R^2)$.

The mapping $U$ is a diffeomorphism of $\overline{B}$ onto
$\overline{D}$ \emph{if and only} if
 \begin{equation}\label{d>0}
 \begin{array}{lll}
 \det DU >0 &\hbox{everywhere on}&\partial B.
 \end{array}
 \end{equation}
 \end{theorem}
 \begin{remark}\label{rem:nonconvex}
 In order to compare this statement with Kneser's Theorem, it is worth noticing that, when $\gamma$ is convex,  (\ref{d>0}) is automatically satisfied.
 Indeed we shall prove, see Lemma \ref{hopf.lemma}, that $\det DU>0$ always holds true on the points of $\partial B$ which are mapped through $\Phi$
  on the part of $\gamma$
 which agrees with its convex hull, see also Definition \ref{nc.def}. As a consequence it is possible to refine the statement of Theorem \ref{main.th},
 by requiring (\ref{d>0}) on a suitable proper subset of $\partial B$. This is the content of Theorem \ref{hopf}. Furthermore, it may be worth stressing that (\ref{d>0}) is,
 in fact, a constraint on the boundary mapping $\Phi$ only. Indeed in Theorem \ref{hopfnc},  by means of the Hilbert
 transform, we shall express  the Jacobian bound  $\det DU>0$  on $\partial B$
 as an explicit, although nonlocal, constraint on the
components of $\Phi$.
  \end{remark}
\begin{remark}\label{rem:shear}
In view of a better appreciation of the strength and novelty of
Theorem~\ref{main.th} let us recall  the so--called method
of \emph{shear construction} introduced by Clunie and Sheil--Small
\cite{css}. Until now, this method has been known \cite[\S
3.4]{duren} as the only other general means for construction of
invertible harmonic mappings, besides Kneser's Theorem. In fact,
we shall show that Theorem \ref{main.th}, and the arguments
leading to its proof, enable us to obtain a new and extremely wide
generalization of the shear construction.  We refer the reader to
Theorem \ref{shearnew} and Corollary \ref{shear:cor} in Section
\ref{sec:shear}, where the shear construction of Clunie and
Sheil--Small is reviewed and our new version is demonstrated.
 \end{remark}
With our next result we return to the original issue for
homeomorphisms. Unfortunately, in this case, the characterization
of the parameterizations $\Phi$, which give rise to homeomorphic
harmonic mappings $U$, is less transparent. It involves the
following classical notion.
\begin{definition}\label{local}
Given $P\in \overline{B}$, a mapping  $U\in C(\overline{B};\mathbb
R^2)$ is a \emph{local homeomorphism} at $P$ if there exists a
neighborhood $G$ of $P$ such that $U$ is one--to--one on $G\cap
\overline{B}$.
\end{definition}
\begin{theorem}\label{hom}
Let $\Phi:\partial B\to \gamma \subset \mathbb R^2$ be a homeomorphism onto a simple closed curve $\gamma$.
Let $D$ be the bounded domain such that $\partial D=\gamma$.
Let $U\in W^{1,2}_{\rm loc} (B;\mathbb R^2)\cap C(\overline{B};\mathbb R^2)$ be the solution to (\ref{main}).

The  mapping $U$ is a homeomorphism of $\overline{B}$ onto $\overline{D}$ \emph{if and only if}, for every $P\in \partial B$,
 the mapping $U$ is a local homeomorphism at $P$.
\end{theorem}
\begin{remark}
Let us note that,  on use of the Riemann Mapping Theorem and the
Caratheodory--Osgood Extension Theorem, see for instance
\cite[Theorem 4.9]{palka}, the disk $B$ can be replaced by any
Jordan domain.
This observation applies also to
 Theorem \ref{main.th} . In this case an analogous result could be stated when the disk $B$ is replaced by any simply connected domain $\Omega$, provided the boundary
of $\Omega$ is smooth enough to guarantee that the map $\omega$ mapping conformally $\Omega$ onto $B$, extends to a
$C^1$ diffeomorphism of   $\overline{\Omega}$ onto $\overline{B}$.
\end{remark}
The paper is organized as follows.

In Section \ref{sec:prelim}
we recall two classical
results of global invertibility, Theorems~\ref{mon.th},
\ref{monc1.th}, and a fundamental result by H. Lewy \cite{lewy}, about
invertible harmonic mappings, Theorem \ref{lewy}.

Section \ref{prel:II} collects a sequence of results which are useful for the proofs of
Theorems~\ref{main.th},~\ref{hom}.  In view of
Theorem \ref{monc1.th} on the inversion of  $C^1$ mappings, our
guiding light towards Theorem \ref{main.th} is to obtain that $\det DU>0$ everywhere in $B$.
This is equivalent to show the absence of critical points for any
linear combination $u_{\alpha} = \cos(\alpha)\, u + \sin(\alpha)\,v$
of the components $u , \ v $ of $U$.
This goal will be achieved
through a number of steps. In Proposition \ref{indice} we show
that, assuming (\ref{d>0}), the number $M_{\alpha}$ of critical
points of $u_{\alpha}$, counted with multiplicities, is finite and
independent of $\alpha$. With Proposition \ref{indolo} we express
the number $M =M_{\alpha}$ in terms of the winding number of the holomorphic
function $f$ the real part of which is $u$. We conclude the Section with
Theorem \ref{wn.th} which enables to compute such winding number
in terms of the boundary mapping $\Phi$.

Section \ref{sec:main} contains the proofs of the main Theorems \ref{main.th}, \ref{hom}.

In  Section \ref{section.H} we present Theorems \ref{hopf},
\ref{hopfnc}, the two refinements of Theorem \ref{main.th} which
we already announced in Remark \ref{rem:nonconvex}.

Section \ref{section:parabola} is mainly devoted to a new proof of
Theorem \ref{parabola.th}. It will be obtained through an
adaptation of an explicit example, which can be traced back at least to
J.C. Wood \cite{wood}, namely, the polynomial harmonic mapping
$F(x,y)=(x, x^2-y^2)$. It is easily seen that such a mapping has a
non--convex range. It  shows also that, contrary to what happens
for holomorphic functions, a harmonic mapping may fail to be open,
see Figure 2  on page \pageref{parabola}. From our construction, we
also obtain that, in Theorem \ref{parabola.th}, the boundary
mapping $\Phi$ can be chosen in such a way that there exists a
curve $\eta \subset B$ on which  $\det DU$ vanishes and such that
$U$ changes its orientation across $\eta$.   In
Remark \ref{fail} we also use this construction to show the
considerable tightness of the condition of local homeomorphism at
the boundary appearing in Theorem~\ref{hom}.

In the final Section \ref{sec:shear}, we  first review the
\emph{shear construction} method of Clunie and Sheil--Small. Then
we state and prove our improved version, namely Theorem
\ref{shearnew}. We conclude with Corollary \ref{shear:cor} which
provides a general construction of harmonic univalent mappings
with prescribed dilatation.

\section{Classical foundations}\label{sec:prelim}
In what follows we shall identify, as usual, points $(x,y) \in
\mathbb{R}^2$ with complex numbers $z=x+iy \in \mathbb{C}$. When
needed, we shall use also  polar coordinates
$z=re^{i\theta}$.

 We now
recall some classical fundamental Theorems which we shall use
several times in the paper.
\begin{theorem}[Monodromy]\label{mon.th}
Let $U\in C(\overline{B};\mathbb R^2)$ be such that

\noindent a) $\Phi=U \big|_{\partial B}$ is a homeomorphism of
$\partial B$ onto a simple closed  curve $\gamma$.

\noindent
b) For every $P\in \overline{B}$, $U$ is a local homeomorphism at $P$.

Then $U$ is a global homeomorphism of $\overline{B}$ onto $\overline{D}$, where $D$ is the bounded domain such that $\partial D=\gamma$.
\end{theorem}
\begin{proof}
A proof can be found in \cite[p.175]{kerek}. Another proof might also be obtained via the theory of light and open mappings of  Stoilow \cite{stoilow}. Results of the same nature in any dimension, but of higher sophistication,  are due to Meisters and Olech \cite{meistol} and Weinstein \cite{weinst}.
\end{proof}
A variant that we shall also use is the following.
\begin{theorem}\label{monc1.th}
Let $U\in C^1(\overline{B};\mathbb R^2)$ be such that

\noindent a$^{\prime}$) $\Phi=U\big|_{ \partial B}$ is a sense
preserving $C^1$ diffeomorphism of $\partial B$ onto a simple
closed  curve $\gamma$.

\noindent
b$^{\prime}$) $\det DU(P)>0$, for every  $P\in \overline{B}$.

Then $U$ is a global diffeomorphism of $\overline{B}$ onto $\overline{D}$.
\end{theorem}
\begin{proof} A proof can be readily obtained as a consequence of Theorem \ref{mon.th}.
\end{proof}
Of a different character is the following Theorem due to H. Lewy
ensuring that harmonic homeomorphisms  are, in fact,
diffeomorphisms as in the holomorphic case.
\begin{theorem}[H. Lewy]\label{lewy}
Let $U:B\to \mathbb R^2$ be harmonic. If $U$ is a sense preserving
homeomorphism, then
\begin{equation*}
\det DU>0 \quad\hbox{everywhere in $B$.}
\end{equation*}
\end{theorem}
\begin{proof}
We refer to \cite{lewy} for a proof.
\end{proof}

\section{Preliminary results}\label{prel:II}
Here we collect some (new) results of essentially topological
nature regarding harmonic functions and harmonic mappings.
\begin{definition}\label{M}
Given a nonconstant harmonic function $u$ defined in $B$, we
denote by $M$ the sum of the multiplicities of its critical
points. Hence $M$ is either a nonnegative integer or $+\infty$.
Given $U=(u,v):B\to \mathbb R^2$ harmonic, we set, for every
 $\alpha\in [0,2 \pi]$
\begin{equation}
u_{\alpha} = \cos(\alpha) u + \sin(\alpha) v
\end{equation}
and denote by $M_{\alpha}$ the sum of the multiplicities of the critical points of $u_{\alpha}$. Our convention is that $M:=M_0$.
\end{definition}
\begin{proposition}\label{indice}
Let $U\in C^1(\overline{B};\mathbb R^2)$ be harmonic in $B$. If $\det DU>0$ on $\partial B$, then for every  $\alpha\in [0,2 \pi]$, the number $M_{\alpha}$ is finite and we have $M_{\alpha}=M$ for every $\alpha \in [0,2 \pi]$.
\end{proposition}
\begin{corollary}\label{ind.cor}
Let $U$ be as in Proposition \ref{indice}. We have $\det DU>0$ everywhere in $B$ if and only if there exists  $\alpha\in [0,2 \pi]$, such that $\nabla u_{\alpha}\neq 0$ everywhere in $B$.
\end{corollary}
\begin{proof}[Proof of Proposition \ref{indice}] Obviously $\nabla u_{\alpha}\neq 0$ everywhere  on $\partial B$ for every $\alpha\in [0,2 \pi]$. By the argument principle for holomorphic functions
\begin{equation*}
M_{\alpha}=\frac{1}{2\pi} \int_{\partial B} \textrm{d} ~{\rm arg}
(\nabla u_{\alpha})\,,\quad \hbox{for every $\alpha\in [0,2
\pi]$.}
\end{equation*}
We shall show that $M_{\alpha}=M_0$ for every $\alpha\in [0,2 \pi]$. It is clear that it suffices to consider
$\alpha\in \left(0,\pi\right)$.
We set
\begin{equation*}J=
\left(
\begin{array}{cc}
0&-1\\
1&0
\end{array}
\right)
\end{equation*}
and we have
\begin{equation*}
\nabla u_{\alpha}\cdot J \nabla u = \sin(\alpha) \det DU>0, \quad \hbox{on $\partial B$,}
\end{equation*}
hence  $|{\rm arg} (\nabla u_{\alpha})-{\rm arg} (J\nabla u)|<\pi$. We conclude that
\begin{equation*}
M_{\alpha}=\frac{1}{2\pi} \int_{\partial B} {\rm d\,arg}(\nabla u_{\alpha})=\frac{1}{2\pi} \int_{\partial B} {\rm d\,arg} (J\nabla u)=M_{0}.
\end{equation*}
\end{proof}
\begin{proof}[Proof of Corollary \ref{ind.cor}] Let us assume that for a given $\alpha\in [0,2 \pi]$, we have  $M_{\alpha}=0$. By Proposition \ref{indice} one has $M_{\alpha}=0$ for every $\alpha\in [0,2 \pi]$. Hence, for every $P\in B$, the vectors $\nabla u(P)$ and $\nabla v(P)$ are linearly independent, that is $\det DU(P)\neq 0$. Being $\det DU>0$ on $\partial B$, by continuity we have $\det DU>0$ everywhere in $B$. The reverse implication is trivial.\end{proof}
\begin{definition} Given a
closed curve $\gamma$, parameterized by $\Phi \in C^1(\partial B;\mathbb R^2)$ and such that
\begin{equation*}
\frac{\partial \Phi}{\partial \theta} \neq 0,\  \textrm{for every} \  \theta\in  [0,2 \pi],
\end{equation*}
we define the {\em winding number} of $\gamma$ as the following integer
\begin{equation*}
{\rm WN}(\gamma) =\frac{1}{2\pi} \int_{\partial B} {\rm d\,arg} \left(\frac{\partial \Phi}{\partial \theta}\right).
\end{equation*}
\end{definition}
\begin{definition}\label{sf}
Let $u$ be a harmonic function in $B$. We denote by $\tilde{u}$ its conjugate harmonic function and we set
\begin{equation*}
f= u + i \tilde{u}.
\end{equation*}
\end{definition}
Note that if, in addition, $u\in C^1(\overline{B})$ and $\nabla u \neq 0$ on $\partial B$, then $f\big|_{\partial B}$ gives us  a regular $C^1$ parametrization of a  closed curve.
\begin{proposition}\label{indolo}
Let $u\in C^1(\overline{B})$ be harmonic in $B$. If $\nabla u\neq 0$ on $\partial B$, then
\begin{equation*}
M={\rm WN}(f(\partial B))-1,
\end{equation*}
with  $M$ as in Definition \ref{M}.
\end{proposition}
\begin{proof} The proof is elementary, and we claim no novelty in this case. We have
\begin{equation*}
{\rm WN}(f(\partial B))=\frac{1}{2\pi} \int_{\partial B} {\rm
d\,arg} \left(\frac{\partial f}{\partial z} \frac{\partial
z}{\partial \theta}\right)= \frac{1}{2\pi} \int_{\partial B}
\textrm{d}~\left[{\rm arg} \left(\frac{\partial f}{\partial z}\right)+
\theta\right]=
\end{equation*}
\begin{equation*}
\frac{1}{2\pi} \int_{\partial B} \textrm{d}~ {\rm arg}(\nabla
u)+1= M+1.
\end{equation*}
\end{proof}
\begin{remark}
If $U=(u,v)\in C^1(\overline{B};\mathbb R^2)$ is such that $\det DU>0$ on $\partial B$, then, for any $P\in \partial B$, the mapping $U$ is a diffeomorphism near $P$. Hence, on $\partial B$,  partial derivatives with respect to $u$ and $v$  make sense.
\end{remark}
\begin{lemma}\label{Jac}
Let $U=(u,v)\in C^1(\overline{B};\mathbb R^2)$ be harmonic in $B$.
If
\begin{equation*}
\det DU>0, \qquad {\rm  on}  \qquad\partial B,
\end{equation*}
then
\begin{equation*}
\frac{\partial \tilde{u}}{\partial v}>0,\qquad  {\rm  on} \qquad \partial B,
\end{equation*}
where   $\tilde{u}$ is the harmonic conjugate of $u$.
\end{lemma}
\begin{proof} We compute
\begin{equation*}
\frac{\partial \tilde{u}}{\partial v} =
\frac{\partial \tilde{u}}{\partial x} \frac{\partial x}{\partial v}+
 \frac{\partial \tilde{u}}{\partial y}\frac{\partial y}{\partial v}=
 \frac{1}{\det DU}
 \left(-\frac{\partial \tilde{u}}{\partial x} \frac{\partial u}{\partial y}+
 \frac{\partial \tilde{u}}{\partial y}\frac{\partial u}{\partial x}
 \right)=
 \frac{|\nabla u|^2}{\det DU}>0.
\end{equation*}
\end{proof}
We are now ready to state a Theorem which contains the main elements towards a proof of Theorem \ref{main.th}.
\begin{theorem}\label{wn.th}
Let $U \in C^1(\overline{B};\mathbb R^2)$ be harmonic in $B$ and let $\Phi=U\big|_{\partial B}$.
If $\det DU>0$ on $\partial B$, then we have
\begin{equation}\label{star}
{\rm WN}( f(\partial B))= {\rm WN} (\Phi(\partial B)).
\end{equation}
\end{theorem}
The proof of Theorem \ref{wn.th} will be based on the following two results.
\begin{proposition}\label{arco}
Given a $C^1$ curve parameterized by $\Phi=(\phi,\psi):[a,b]\to \mathbb R^2$ and such that $\phi^{\prime}\neq 0$ in $(a,b)$ and $\phi^{\prime}(a)=\phi^{\prime}(b)= 0$ and given a
$C^1$ function $g:\psi([a,b])\to \mathbb R$ with $g^{\prime}>0$ in $\psi((a,b))$, consider the curve $\widetilde{\Phi}:[a,b]\to \mathbb R^2$ given by $\widetilde{\Phi}=(\phi,g(\psi))$.
We have
\begin{equation}\label{arg}
\int_a^b {\rm d\,arg}(\widetilde{\Phi}^{\prime})=
\int_a^b {\rm d\,arg}(\Phi^{\prime}).
\end{equation}
\end{proposition}
\begin{lemma}\label{wn.prop}
Under the same assumptions of Theorem \ref{wn.th}, assuming in addition that $\frac{\partial u}{\partial \theta}\Big|_{\partial B}$ vanishes at finitely many points, we have that (\ref{star}) holds.
\end{lemma}
\begin{proof}[Proof of Proposition \ref{arco}]
Without loss of generality we may assume $\phi^{\prime}>0$ in $(a,b)$. We have that
both ${\rm arg} \big(\Phi^{\prime}\big)$ and ${\rm arg} \big(\widetilde{\Phi}^{\prime}\big)$ take values in $(-\pi,\pi)$. Hence
\begin{equation*}
\int_a^b {\rm d\,arg}(\Phi^{\prime})={\rm arg}\big(\Phi^{\prime}(b^-)\big)-{\rm arg}\big(\Phi^{\prime}(a^+)\big)
\end{equation*}
and also
\begin{equation*}
\int_a^b {\rm d\,arg}(\widetilde{\Phi}^{\prime})={\rm
arg}\big(\widetilde{\Phi}^{\prime}(b^-)\big)-{\rm
arg}\big(\widetilde{\Phi}^{\prime}(a^+)\big).
\end{equation*}
Now, we compute
\begin{equation*}
{\rm arg}(\Phi^{\prime}(b^-))= {\rm arg}(\widetilde{\Phi}^{\prime}(b^-))=\pm\frac{\pi}{2},
\end{equation*}
and also
\begin{equation*}
{\rm arg}(\Phi^{\prime}(a^+))= {\rm arg}(\widetilde{\Phi}^{\prime}(a^+))=\pm\frac{\pi}{2}.
\end{equation*}
Hence (\ref{arg}) follows.\end{proof}
\begin{proof}[Proof of Lemma \ref{wn.prop}]
Up to a rotation in the $x,y$ coordinates, we may assume without loss of generality that there exists a partition of $[0,2\pi]$, $0=\theta_0<\theta_1<\hdots<\theta_N=2\pi$ such that
\begin{equation*}
\frac{\partial \phi}{\partial \theta_k} (\theta)  = 0  \quad\hbox{for all}\quad k=0,1,\hdots,N-1,
\end{equation*}
and
\begin{equation*}
\frac{\partial \phi}{\partial \theta} (\theta)\neq 0 \quad\hbox{in}\quad (\theta_k,\theta_{k+1}),\quad\hbox{for every}\quad k=0,1,\hdots,N-1.
\end{equation*}
On each interval $[\theta_k,\theta_{k+1}]$, we have
\begin{equation*}
f(e^{i \theta})=\Big(\phi(\theta),g(\psi(\theta))\Big) \quad\hbox{with}\quad g\big(\psi(\theta)\big)=\tilde{u}(e^{i \theta})
\end{equation*}
and, by Lemma \ref{Jac},
\begin{equation*}
\frac{\partial g}{\partial \psi}=\frac{\partial \tilde{u}}{\partial v}>0.
\end{equation*}
Hence, by Proposition \ref{arco}
\begin{equation*}
\int_{\theta_k}^{\theta_{k+1}} {\rm d\,arg} \left(\frac{\partial f}{\partial \theta}\right)=
\int_{\theta_k}^{\theta_{k+1}} {\rm d\,arg} \left(\frac{\partial \Phi}{\partial \theta}\right)\quad\hbox{for every}\quad k=0,1,\hdots,N-1
\end{equation*}
and (\ref{star}) follows.\end{proof}
\begin{proof}[Proof of Theorem \ref{wn.th}]
By continuity, there exists $\rho\in (0,1)$, such that $\det DU
>0$ in $\overline{B}\setminus B_{\rho}(0)$ and, consequently,
$\frac{\partial f}{\partial z}\neq 0$ in $\overline{B}\setminus
B_{\rho}(0)$. Therefore the numbers
\begin{equation*}
{\rm WN}\Big(f\big(\partial B_r(0)\big)\Big) \quad\hbox{and}\quad{\rm WN}\Big(U\big(\partial B_r(0)\big)\Big)
\end{equation*}
are constant with respect to $r\in[\rho,1]$. Since $
u\big|_{\partial B_r}(\theta)$ is a nonconstant real analytic
function of $\theta$, we have that  $\frac{\partial u}{\partial r}
\big(r e^{i \theta}\big)$ vanishes at most on a finite set of
angles $\theta_j \in [0,2 \pi]$. Applying Lemma \ref{wn.prop}  to
$U(r\cdot)$ rather than $U$, we obtain
\begin{equation*}
{\rm WN}\Big(f\big(\partial B\big)\Big) =
{\rm WN}\Big(f\big(\partial B_r(0)\big)\Big) =
{\rm WN}\Big(U\big(\partial B_r(0)\big)\Big)=
{\rm WN}\Big(\Phi\big(\partial B\big)\Big).
\end{equation*}
\end{proof}

\section{Proofs of the main Theorems}\label{sec:main}
\begin{proof}[Proof of Theorem \ref{main.th}]
Let us assume that (\ref{d>0}) holds. By assumption $\Phi$ is
one--to--one and sense preserving. Hence by Theorem \ref{wn.th},
\begin{equation*}
{\rm WN}\Big(f\big(\partial B\big)\Big) =
{\rm WN}\Big(\Phi\big(\partial B\big)\Big)=1.
\end{equation*}
By Proposition \ref{indolo}, $\nabla u$ never vanishes in $B$.  By
Corollary \ref{ind.cor},  $\det DU>0$ everywhere in
$\overline{B}$. By Theorem \ref{monc1.th}, $U:\overline{B}\to
\overline{D}$ is a diffeomorphism.
The reverse implication is obvious.
\end{proof}
Our next goal is to prove Theorem \ref{hom}. We need the following preliminary Lemma.
\begin{lemma}\label{top}
Assume $\Phi:\partial B\to \gamma\subset \mathbb R^2$ is a
homeomorphism onto a simple closed curve $\gamma$. Let $U\in C^2
(B;\mathbb R^2)\cap C(\overline{B};\mathbb R^2)$ be  the
solution to   (\ref{main}). If, in addition, for every $P\in
\partial B$ the mapping $U$ is a local homeomorphism near $P$,
then there exists $\rho\in (0,1)$ such that $U$ is a
diffeomorphism of $B\setminus \overline{B_{\rho}(0)}$ onto
$U\Big(B\setminus \overline{B_{\rho}(0)}\Big)$.
\end{lemma}
\begin{proof}
By the compactness of $B$, there exist finitely many points $P_1,\hdots,P_k\in \partial B$ and a number $\delta>0$ such that
\begin{equation*}
\partial B \subset  \bigcup\limits_{k=1}^K B_{\delta}(P_k),
\end{equation*}
and $U$ is one--to--one on $B_{2 \delta}(P_k)\cap\overline{B}$
for every $k$. Note that there exists $\rho\in (0,1)$ such that
\begin{equation*}
\overline{B}\setminus B_{\rho}(0)\subset \bigcup\limits_{k=1}^K
B_{\delta}(P_k).
\end{equation*}
Let $P,Q$ be two distinct points in $\overline{B}\setminus
B_{\rho}(0)$. If $|P-Q|<\delta$, then there exists $k=1,\hdots,K$
such that $P,Q\in B_{2 \delta}(P_k)$ and, hence, $U(P)\neq U(Q)$.
Assume now $|P-Q|\geq \delta$. Let
\begin{equation*}
P^{\prime}=\frac{P}{|P|}\quad,\quad Q^{\prime}=\frac{Q}{|Q|}.
\end{equation*}
We have $|P-P^{\prime}|<1-\rho,\, |Q-Q^{\prime}|<1-\rho,$ and thus
\begin{equation*}
|P^{\prime}- Q^{\prime}| >|P-Q| - 2(1-\rho) \geq \delta -2(1-\rho).
\end{equation*}
Choosing $\rho$ such that $(1-\rho)<\frac{\delta}{4}$, we have
$|P^{\prime}- Q^{\prime}|>\frac{\delta}{2}.$ Now we use the fact
that $P^{\prime}$ and $Q^{\prime}$ belong to $\partial B$ and
$\Phi$ is one--to--one to deduce that there exists $c>0$ such that
\begin{equation*}
|\Phi(P^{\prime})- \Phi(Q^{\prime})|   \geq c.
\end{equation*}
Recall that $U$ is uniformly continuous on $\overline{B}$. Denoting by $\omega$ its modulus of continuity, we have
\begin{equation*}
|U(P) - U(Q) |\geq |U(P^{\prime})- U(Q^{\prime})| -2 \omega(1-\rho)=
\end{equation*}
\begin{equation*}
|\Phi(P^{\prime})- \Phi(Q^{\prime})| -2 \omega(1-\rho)\geq c -2 \omega(1-\rho).
\end{equation*}
Choosing $\rho$ such that $1-\rho<\omega^{-1}\big(\frac c 4\big)$ we obtain
\begin{equation*}
|U(P) - U(Q) |\geq \frac c 2 >0,
\end{equation*}
which implies the injectivity of $U$ in $\overline{B}\setminus
B_{\rho}(0)$. Consequently, by Theorem~\ref{lewy}, $\det DU\neq 0$
in $B\setminus \overline{B_{\rho}(0)}$ and the thesis follows.
\end{proof}
\begin{proof}[Proof of Theorem \ref{hom}]
We assume that, for every $P\in \partial B$, $U$ is a local homeomorphism near $P$ and prove that $U:\overline{B}\to \overline{D}$ is a homeomorphism. The opposite implication is trivial. In view of Theorem \ref{mon.th} it suffices to show that $\det DU\neq 0$ everywhere in $B$.

For every $r\in (0,1)$, let us write $\Phi^r:\partial B\to \mathbb R^2$ to denote the application given by
\begin{equation*}
\Phi^r(e^{i\theta})= U(r e^{i\theta}),\quad \theta\in [0,2\pi].
\end{equation*}
By Lemma \ref{top}, there exists $\rho\in (0,1)$ such that for every $r\in (\rho,1)$ the mapping $\Phi^r:\partial B\to \gamma_r\subset \mathbb R^2$
is a diffeomorphism of $\partial B$ onto a simple closed curve $\gamma_r$, and $U(r\cdot)\in C^1(\overline{B};\mathbb R^2)$
solves (\ref{main}) with $\Phi$ replaced by $\Phi^r$.
Then, by Theorem \ref{main.th} we obtain
\begin{equation*}
\det DU(r z)\neq 0, \quad\hbox{for every}\quad z\in \overline{B}
\end{equation*}
that is
\begin{equation*}
\det DU(z)\neq 0, \quad\hbox{for every}\quad z\in \overline{B_r(0)}.
\end{equation*}
Finally, by Lemma \ref{top} we have $\det DU\neq 0$ in $B\setminus
\overline{B_{\rho}(0)}$ so that $\det DU\neq 0$ everywhere in $B$.
\end{proof}

\section{Variations upon Theorem \ref{main.th}}\label{section.H}
Let us introduce some definitions borrowed from the literature on
minimal surfaces \cite{sauvigny}.
\begin{definition}\label{nc.def}
Given a Jordan domain $D$, let us denote by ${\rm co} (D)$ its
convex hull. We define the {\em convex part} of $\partial D$ as
the closed set $\gamma_c=\partial D\cap \partial ({\rm co}(D))$.
Consequently we define  the {\em non--convex part} of $\partial D$
as the open set $\gamma_{nc}=\partial D\setminus \partial ({\rm
co}(D))$.
\end{definition}
\begin{theorem}\label{hopf}
Under the same assumptions as in Theorem \ref{main.th}, the
mapping $U$ is a diffeomorphism of $\overline{B}$ onto
$\overline{D}$ \emph{if and only if}
\begin{equation}\label{ncd>0}
\det DU>0\quad \hbox{everywhere on $\Phi^{-1}(\gamma_{nc})$},
\end{equation}
where $\gamma_{nc}$ is the set introduced in Definition
\ref{nc.def} above.
\end{theorem}
First we prove the following Lemma.
\begin{lemma}\label{hopf.lemma}
Under the  assumptions of Theorem \ref{main.th}, we always have
\begin{equation}\label{d>0nc.eq}
\det DU>0\quad\hbox{everywhere on $\Phi^{-1}(\gamma_{c})$}.
\end{equation}
\end{lemma}
\begin{proof}
Let $P\in \Phi^{-1}(\gamma_c)$ and $Q=\Phi(P)$. Let $l$ be a
support line for ${\rm co}(D)$ at $Q$.
Without loss of generality, we may assume
\begin{equation*}
\begin{array}{lll}
l=\{(u,v)\in \mathbb R^2: v=0\}&,&{\rm co}(D)\subset \{(u,v)\in
\mathbb R^2: v>0\}.
\end{array}
\end{equation*}
Thus the second component $\psi$ of $\Phi$ satisfies
\begin{equation*}
\psi\geq 0 \quad \hbox{everywhere},\quad \psi(P)=0,
\end{equation*}
hence $P$ is a minimum point for $\psi$ and therefore
\begin{equation*}
\frac{\partial \psi}{\partial \theta}(P)=0 .
\end{equation*}
Moreover, being $\Phi$ orientation preserving, we have that $\phi$
is increasing at $P$ and also
\begin{equation*}
\left(\frac{\partial \phi}{\partial
\theta}(P)\right)^2=\left|\frac{\partial \Phi}{\partial
\theta}(P)\right|^2>0
\end{equation*}
so that
\begin{equation*}
\frac{\partial \phi}{\partial \theta}(P)>0.
\end{equation*}
On the other hand,  Hopf's Lemma gives
\begin{equation*}
\frac{\partial v}{\partial r}(P)<0.
\end{equation*}
Consequently
\begin{equation*}
\det DU(P)= -\frac{\partial \phi}{\partial
\theta}(P)\frac{\partial v}{\partial r}(P)>0.
\end{equation*}
\end{proof}
\begin{proof}[Proof of Theorem \ref{hopf}]
The proof is a straightforward consequence of Theorem
\ref{main.th} and of the above Lemma \ref{hopf.lemma}.
\end{proof}
We now turn to the Hilbert transform formalism. For any $g\in
L^2((0,2\pi))$, let
\begin{equation}\label{ht}
\mathcal{H}g(\theta):= \frac{1}{2 \pi}\ {\rm P.V.}\!\!\! \int_0^{2
\pi} \frac{g(\tau)}
 {
 {\rm tan}
  \left(
  \frac{\theta -\tau}{2}
  \right)
  } d\tau,\qquad\theta\in[0,2\pi],
\end{equation}
be the Hilbert transform on the unit circle, see for instance
\cite[p. 145]{smirnov}. The following is an equivalent formulation of
Theorem \ref{hopf} and thus of Theorem \ref{main.th}.
\begin{theorem}\label{hopfnc}
Under the same assumptions as in Theorem \ref{main.th}, $U$ is a
diffeomorphism of $\overline{B}$ onto $\overline{D}$ \emph{if and
only if} the components $\phi$ and $\psi$ of $\Phi$ satisfy
\begin{equation}\label{htd>0}\frac{\partial \phi}{\partial \theta}  \mathcal{H}\left(\frac{\partial
\psi}{\partial \theta}\right) -\frac{\partial \psi}{\partial
\theta}  \mathcal{H}\left(\frac{\partial \phi}{\partial
\theta}\right)>0 \quad\hbox{everywhere on
$\Phi^{-1}(\gamma_{nc})$.}
\end{equation}
\end{theorem}
\begin{proof} Expressing $\det DU$ in polar coordinates we have, on $\partial B$,
\begin{equation*}
\det DU= \frac{\partial u}{\partial r} \frac{\partial v}{\partial
\theta}- \frac{\partial u}{\partial \theta}\frac{\partial
v}{\partial r}.
\end{equation*}
Since $\frac{1}{r} \frac{\partial u}{\partial \theta}$ is the
harmonic conjugate of $\frac{\partial u}{\partial r}$, we have
that
\begin{equation}
\frac{\partial u}{\partial r} =-\mathcal{H} \left(\frac{\partial
u}{\partial \theta}\right)  \quad\hbox{everywhere on $\partial B$}
\end{equation}
and the same formula, obviously applies for $v$. By assumption
$u,v\in C^1(\overline{B})$, hence we obtain
\begin{equation*}
\det DU= \frac{\partial \phi}{\partial \theta}
\mathcal{H}\left(\frac{\partial \psi}{\partial \theta}\right)
-\frac{\partial \psi}{\partial \theta}
\mathcal{H}
\left(\frac{\partial \phi}{\partial
\theta}\right)\quad\hbox{everywhere on $\partial B$}.
\end{equation*}
Hence condition (\ref{ncd>0}) is equivalent to (\ref{htd>0}).
\end{proof}

\section{The counterexample}\label{section:parabola}
\begin{proof}[Proof of Theorem \ref{parabola.th}] It suffices to prove the Theorem with $D$ replaced by $T D$ where $T$ is an invertible
affine transformation. In fact, the Theorem will be proved with $\Phi$ and $U$ replaced by $T^{-1}\Phi$ and $T^{-1}U$
respectively.

If $D$ is not convex, we can find a support line $l$ of its convex
hull ${\rm co} (D)$ which touches $\partial D$  on (at least)  two
points $A$ and $B$ and such that the {\em open} segment
$\overline{AB}$ is outside $\overline{D}$. The midpoint $C$ of
$\overline{AB}$ is at a positive distance from $\overline{D}$. We
can also find $E\in {\rm co} (D)\setminus
 D$ such that the segment $\overline{CE}$ is
perpendicular to $\overline{AB}$ and it lies outside
$\overline{D}$.

Next we consider $K$, the largest closed cone, with vertex at $E$,
such that $K\cap D =\emptyset$. Note that $K\subseteq K^{\prime}$,
where $K^{\prime}$ is the convex cone with vertex at $E$ and such
that
 $A,B\in \partial K^{\prime}$. Therefore the cone $K$ is convex.
Let $\alpha,\beta$ be the half--lines such that $\alpha\cup \beta
= \partial K$. Then $\alpha$ intersects $\partial D$ in at least
one point $A^{\prime}$ and similarly $\beta$ intersects $\partial
D$ in at least one point $B^{\prime}$. Up to an affine
transformation, we may assume that
$|A^{\prime}-E|=|B^{\prime}-E|$.

Let ${\mathcal P}$ be the unique parabola contained in $K$ which passes through $A^{\prime}$ and
$B^{\prime}$. Up to a further affine transformation, we may assume
\begin{equation*}
{\mathcal P}=\{(u,v)\big| v=u^2\},\quad A^{\prime}=(p,p^2),\quad B^{\prime}=(-p,p^2)\quad\hbox{for some}\quad p>0.
\end{equation*}
Consider the harmonic mapping $F:\mathbb R^2\to \mathbb R^2$ given by
\begin{equation*}
F(x,y)= (x,x^2-y^2).
\end{equation*}
Set
\begin{equation*}
Y_+=\{y\geq 0\},
\quad Y_-=\{y\leq 0\}
\quad \hbox{and}\quad V_-=\{(u,v)\in \mathbb R^2\big| v\leq u^2\}.
\end{equation*}
The mappings $F_{\pm}=F\Big|_{Y_{\pm}}:Y_{\pm}\to V^-$ are both
one--to--one.

Let $\gamma_1,\gamma_2$ be the simple open arcs of
$\gamma=\partial D$ whose endpoints are $A^{\prime}$ and
$B^{\prime}$. Consider the closed curve $\Gamma$ obtained by
gluing together the arcs
\begin{equation*}
\Gamma_1=F^{-1}_+(\gamma_1),\quad\Gamma_2=F^{-1}_-(\gamma_2),
\end{equation*}
through their {\em common} endpoints  $F^{-1}(A^{\prime})$ and $F^{-1}(B^{\prime})$.
Then $\Gamma$ is a simple closed curve which intersects the line $\{y=0\}$  exactly at the two points,
 $F^{-1}(A^{\prime})=(p,0)$ and $F^{-1}(B^{\prime})=(-p,0)$.
Let $G$ be the Jordan domain bounded by $\Gamma$ and let $\omega$
be a conformal mapping  $\omega :B \to G$ which extends to a
homeomorphism of $\overline{B}$ onto $\overline{G}$. We define
\begin{equation}\label{Uexample}
\Phi= (F \circ \omega)\big|_{\partial B},\quad\hbox{and}\quad U=F\circ \omega,\quad\hbox{in}\quad B.
\end{equation}
One then verifies that $\Phi:\partial B\to \partial D$ is a
homeomorphism, that $U$ solves (\ref{main}) and that it is not
one--to--one. In fact, $\det DU$ changes its sign across the curve
\begin{equation*}
\eta=\omega^{-1}\Big(\big\{(x,0)\in \mathbb R^2 : |x|<p\big\}\Big).
\end{equation*}
Moreover, $\det DU=0$ in $B$ if and only if $(x,y)\in \eta$ and
$U$ maps the curve $\eta$    in a one--to--one way onto  the arc
of the parabola ${\mathcal P}$ which joins $A^{\prime}$ to
$B^{\prime}$ and which lies outside $\overline{D}$.
\end{proof}
\begin{figure*}[h]
\label{construction}
\centerline{
\includegraphics[scale=0.3,angle=0]{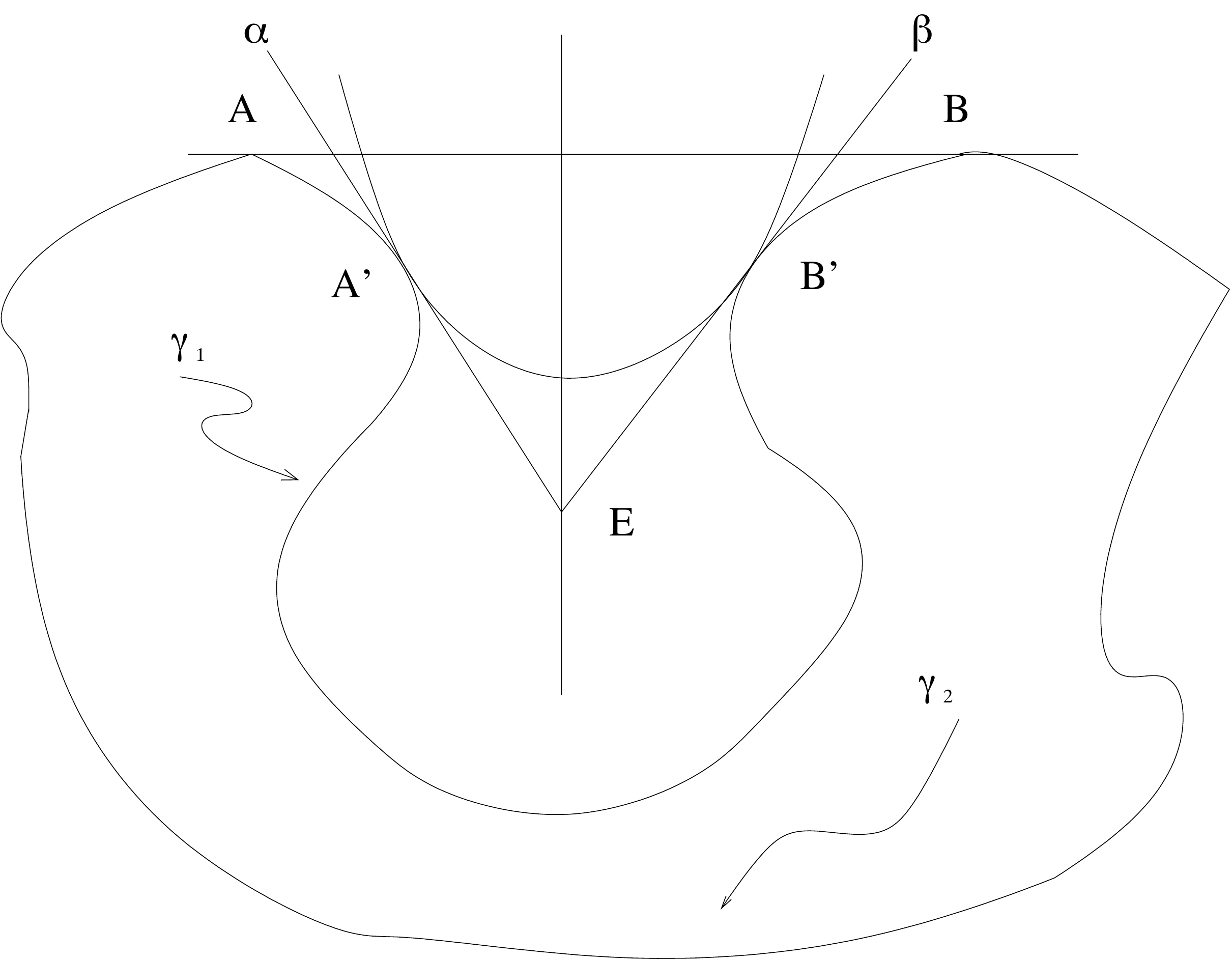}
}
\caption{A sketch of the construction of the counterexample}
\end{figure*}
\begin{figure}[h]
\label{parabola}
\centerline{
\includegraphics[scale=0.3,angle=0]{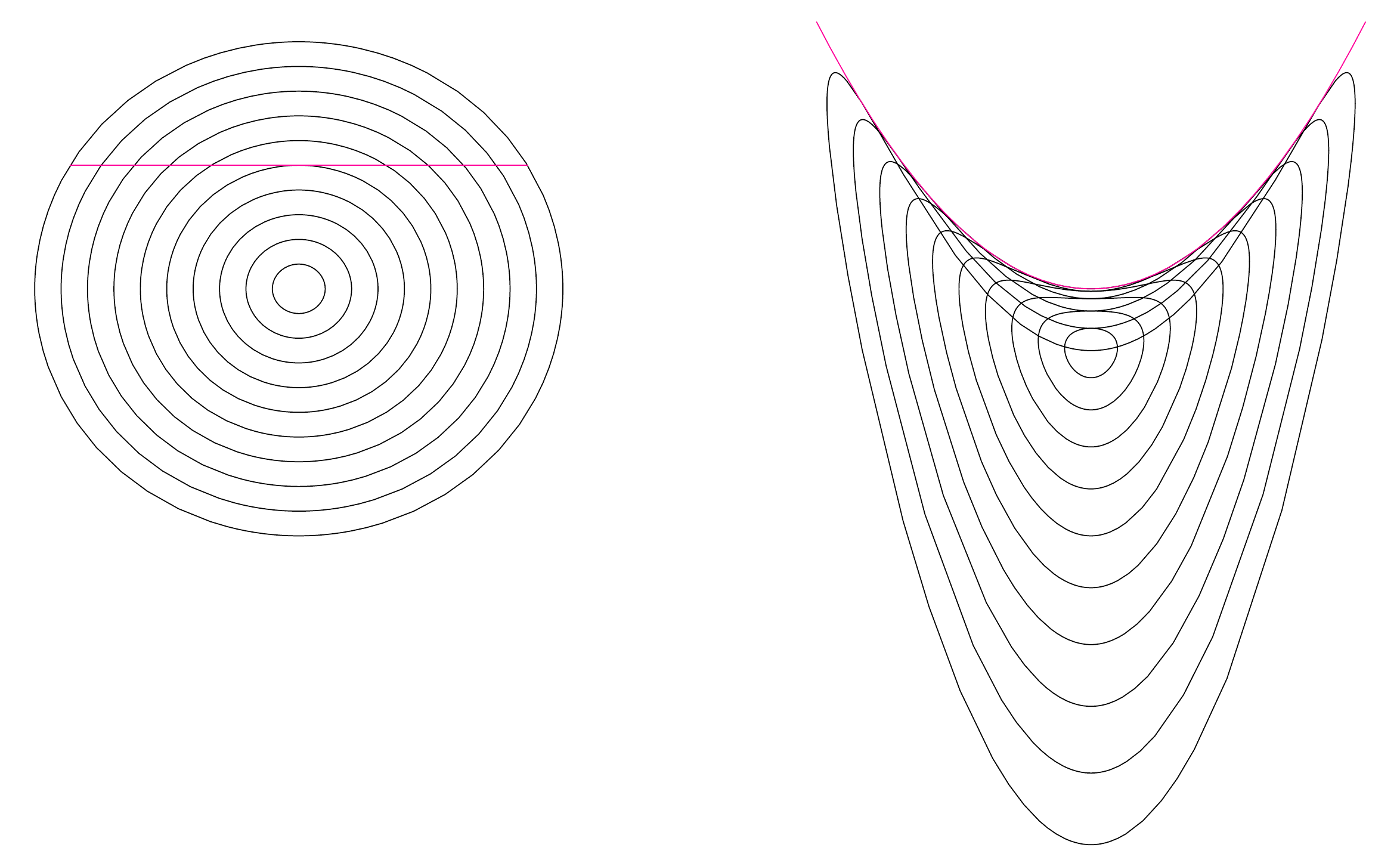}
}
\caption{The basic example: the mapping $F(x,y)=(x,x^2-y^2)$}
\end{figure}
\begin{remark}\label{fail}
In the construction of our counterexample the harmonic mapping $U$ given by (\ref{Uexample}) fails to be  a local homeomorphism on $\partial B$ \emph{exactly}  at the points $A^{\prime\prime}=\omega^{-1}(p,0)$, $B^{\prime\prime}=\omega^{-1}(-p,0)$.
This is a clear indication of how close to optimal
Theorem \ref{hom} is. In fact, the conclusion of Theorem \ref{hom} does not hold if  the condition
 \begin{quotation}
 \emph{for every $P\in \partial B$,
 the mapping $U$ is a local homeomorphism at $P$}
 \end{quotation}
 is relaxed to
\begin{quotation}
 \emph{for every $P\in \partial B$, \emph{except} possibly at two points,
 the mapping $U$ is a local homeomorphism at $P$.}
 \end{quotation}
\end{remark}

\section{The shear construction revisited} \label{sec:shear}

Let us recall the so--called   \emph{shear
construction} method due to Clunie and Sheil-Small \cite{css}. In order
to conform to the language of the previous sections, we shall
adapt their definitions to the current notation of this paper.

Let $U=(u,v)$ be a harmonic mapping on $B$ and let $\tilde u$ and
$\tilde v$ be the harmonic conjugates of $u$ and $v$ respectively.
We already introduced the holomorphic function
\begin{equation}\label{folomorfa}
f=u+i \tilde u \ ,
\end{equation}
accordingly, we define
\begin{equation}
g=v+i\tilde v \ .
\end{equation}
Let us further introduce the following linear combinations
\begin{equation}
G=\frac 1 2\left(f+ig\right) \ ,\ H=\frac 1 2 \left(f-ig\right) .
\end{equation}
Then we have
\begin{equation}\label{can-rep}
U=G+\overline{H}
\end{equation}
which is usually  called the \emph{canonical representation} of $U$. Note that, by construction,
we have that $f$ as defined by (\ref{folomorfa}), satisfies
\begin{equation}\label{can-rep-f}
f=G+H\ .
\end{equation}
Here
with slight, although customary, abuse of notation we have
identified
 $U =(u,v)$ with $u + i v$.
\begin{definition}\label{cvd}
For any $\theta \in [0,\pi)$, a set $K\subseteq
\mathbb R^{2}$ is called {\it convex in the  direction $\zeta=e^{i\theta}$}, if
any line parallel to $\zeta$ intersects $K$ in a connected set,
possibly empty or unbounded. We denote by $C_{\theta}$ the class of
such sets. In particular, $C_{\pi/2}$ denotes the class of  sets which are
convex in the \emph{vertical} direction and we write $f(B)\in C_{\theta}$ to indicate that the
range of $f$ is convex in the direction $e^{i \theta}$.

\end{definition}

The basic Theorem of the shear construction method is as follows.
\begin{theorem}[Clunie and Sheil-Small]\label{shearold}
Let $U$ be a harmonic mapping on $B$ with canonical representation
as in (\ref{can-rep}), let $f$ be defined by (\ref{can-rep-f}) and assume that
\begin{equation}
\det DU > 0\,\,  \mbox{in}\,\, B\ .
\end{equation}
The following two conditions are equivalent
\begin{equation}\label{Ucvd}
U \text{ is one--to--one and } U(B)\in C_{\pi/2}\ ,
\end{equation}
\begin{equation}\label{fcvd}
f\text{ is one--to--one and } f(B)\in C_{\pi/2}\ .
\end{equation}
\end{theorem}
A slightly more involved version of this result is available,
which applies when the class of sets convex in one direction is
replaced with the class of \emph{close--to-convex} sets, we refer
to \cite{css, bh} for the definition and details. Our new version
is the following.

\begin{theorem}\label{shearnew}
Let $U\in C^2 (B;\mathbb R^2)\cap C^1(\overline{B};\mathbb R^2)$
be harmonic on $B$ with canonical representation
as in (\ref{can-rep}), let $f$ be defined by (\ref{can-rep-f}) and assume  that
\begin{equation}
\det DU > 0\,\,  \mbox{on}\,\, \partial B \ ,
\end{equation}
then the following conditions are equivalent
\begin{equation}\label{Ubordo}
U \big|_{\partial B} \text{ is one--to--one }\ ,\end{equation}
\begin{equation}\label{fbordo}
f \big|_{\partial B} \text{ is one--to--one }\ ,
\end{equation}
\begin{equation}\label{Udentro}
U  \text{ is one--to--one on }\overline{B} \ ,
\end{equation}
\begin{equation}\label{fdentro}
f  \text{ is one--to--one on }\overline{B} \ .
\end{equation}
\end{theorem}
\begin{proof}[Proof (sketch).] By Theorem
\ref{main.th}, \eqref{Ubordo} is equivalent to \eqref{Udentro}.
Theorem \ref{main.th} can also be used to show that \eqref{fbordo}
is equivalent to \eqref{fdentro}, however a more direct proof can
be obtained by the use of the classical argument principle. From
Theorem \ref{wn.th} we readily obtain that  \eqref{Ubordo} is
equivalent to \eqref{fbordo}.
\end{proof}
Observe that, in comparison to Theorem \ref{shearold}, at the minor
price of assuming $C^1$ regularity up to the boundary for $U$, we
have obtained the advantage that the condition of non--vanishing
of the  Jacobian  is now required \emph{on the boundary only}, and
that we do not need anymore the assumption of convexity in some
direction.
We recall also that one of the main interest of Theorem
\ref{shearold} is that it allows to construct univalent harmonic
functions with prescribed \emph{dilatation}
\begin{equation}\label{omegaolomorfa}
\omega = \frac{\overline{U_{\bar{z}}}}{U_{z}} \ .
\end{equation}
We refer the reader to the monograph of P.  Duren \cite{duren} for
more details about the meaning of the dilatation of a harmonic
mapping, also called \emph{second complex dilatation} or
\emph{analytic dilatation}. For the present purposes it suffices
to recall that $\omega$ is holomorphic and that, at any point, the
condition $\det DU
> 0$ is equivalent to $|\omega|<1$. If we are given a univalent holomorphic
function $f$ and a holomorphic function $\omega$ such that
$|\omega|<1$ in $B$ {\it and} such that  $f(B)\in C_{\pi/2}$, then
one can construct a harmonic univalent function $U$ such that
$U(B)\in C_{\pi/2}$
 and which has the canonical
representation $U= G+\overline{H}$ where $G$ and $H$ are
determined by the linear system
\begin{equation}\label{system}
\left\{
\begin{array}{lll}
G_{z}+H_{z}&=&f_{z}\\
\omega G_{z}-H_{z}&=&0.
\end{array}
\right.
\end{equation}
In this way one obtains a harmonic injective mapping with
prescribed dilatation $\omega$. The name of \emph{shear
construction} is  related to the mechanical concept of \emph{shear
deformation}. Indeed $U$ is obtained from $f$, by keeping one
component fixed (in this case the real part) and by deforming the
other (in this case the imaginary part). The fundamental drawback
is that one cannot apply the method when no sort of convexity
assumption on the range of $f$ is available. As a consequence of
our new Theorem \ref{shearnew} we can remove this kind of
requirement.

\begin{corollary}\label{shear:cor} Let $f , \omega$ be holomorphic
functions in $B$ such that $f$ extends to a $C^1$ invertible
mapping on $\overline{B}$, $\omega$ extends continuously to
$\overline{B}$ and it satisfies
$$
|\omega|<1 \ , \text{ in }  \overline{B} \ .
$$
Then, given $G, H$ the holomorphic solutions to \eqref{system},
the harmonic mapping $U= G+\overline{H}$ is a diffeomorphism on
$\overline{B}$, it satisfies $\mathfrak{Re} U = \mathfrak{Re} f$
and its dilatation equals $\omega$ in $\overline{B}$.
\end{corollary}
\begin{proof} Straightforward consequence of Theorem
\ref{shearnew}. \end{proof}

\begin{remark}It is evident that if we merely
assume $f , \omega$ be holomorphic functions in $B$ such that $f$
is invertible on $B$ and $\omega$
 satisfies
$$
|\omega|<1 \ , \text{ in }  B \ ,
$$
then the same construction yields a harmonic mapping $U$ which is
a diffeomorphism on the open disk $B$.  In fact, similarly to what
we did in the proof of Theorem \ref{hom}, it suffices to apply
Corollary \ref{shear:cor} by shrinking the independent variable
$z\in B$  to $rz \in B_r(0)$ for any $0<r<1$. It is also evident,
indeed, that the above construction provides a complete
characterization of harmonic diffeomorphisms. In fact, given a
harmonic diffeomorphism $U$ either on $\overline{B}$, or on $B$,
$f$ and $\omega$ are immediately obtained by \eqref{folomorfa},
\eqref{omegaolomorfa}.

\end{remark}

\begin{acknowledgement} The authors express their gratitude to
G. F.  Dell'Antonio for inspiring and fruitful conversations
occurred while completing this paper.
\end{acknowledgement}
\bibliographystyle{amsplain}
\providecommand{\bysame}{\leavevmode\hbox
to3em{\hrulefill}\thinspace}
\providecommand{\MR}{\relax\ifhmode\unskip\space\fi MR }
\providecommand{\MRhref}[2]{%
  \href{http://www.ams.org/mathscinet-getitem?mr=#1}{#2}
} \providecommand{\href}[2]{#2}

\end{document}